\documentstyle[12pt]{article}

\newcommand{\spone}{0.9}
\newcommand{\sphalf}{1.15}
\newcommand{\sptwo}{1.4}

\newcommand{\singlespace}{\edef\baselinestretch{\spone}\Large\normalsize}
\newcommand{\kimspace}{\edef\baselinestretch{\sphalf}\Large\normalsize}
\newcommand{\doublespace}{\edef\baselinestretch{\sptwo}\Large\normalsize}

\doublespace
\begin{document}

\begin{center}
{\bf Odd Order Pandiagonal Latin and Magic %Squares and
Cubes
in Three and Four Dimensions}\\
Solomon Gartenhaus, Purdue University, W. Lafayette, IN 47907\\
garten@physics.purdue.edu
\end{center}
\vspace{7pt}

\noindent
{\bf I  Introduction}
\vspace{5pt}

%We define
By a magic square of order $n$ is here meant an arrangement,
without repeats, of the integers
$\{0,1,2,...n^2$ - 1\} into the $n^2$ cells of an $n \times n$ square
in a way that the sum of the elements of each row, of each column and of each
of the two diagonals is the same.  Since the sum of
$1 + 2 + ... + (n^2-1)$ is
%0, 1, 2,.. n$^2$ - 1 is
$n^2(n^2$-1)/2 and since for a magic square this must be equal to the sum of
the sum of the $n$-integers in each of the $n$-rows (or columns) it follows that the
common sum, $\sigma_2$ of these integers
must be
$$ \sigma_2 = \frac{n(n^2 - 1)}{2}
\eqno{(1.1)}
$$

Figure I shows two magic squares for $n$ = 4 and $n$ = 5
\vspace{1pt}

$$
\begin{array}{rllrllllr}
0 & 10 & 15 & 5~~~~~~~~~~~~~~&13 & 19 & 20 & 1 & 7 \\
7 & 13 & 8  & 2~~~~~~~~~~~~~~&21 & 2 & 8 & 14 & 15 \\
9 &  3 & 6 & 12 ~~~~~~~~~~~~~~&9 & 10 & 16& 22 & 3 \\
14 & 4 & 1 & 11~~~~~~~~~~~~~~& 17 & 23 & 4 & 5 & 11 \\
& & &    & 0 & 6 & 12 & 18 & 24 \\
& (a) & &~~~~~~~~~~~~~~& & & (b) & \\
\end{array}
$$
\vspace{1pt}

\begin{center}
Figure I
\end{center}
\vspace{3pt}

\noindent
with sum values 30 and 60, respectively.
Of particular interest is the magic square
in Ib.  Not only do the sums along the rows, the columns and the two diagonals have
the value 60, but so do the eight {\bf broken diagonals}
that are obtained if all partial diagonals are completed by imagining the
square bent around into a cylinder.  Examples of broken diagonals in
Figure Ib are:
21, 10, 4, 18, 7; 9,
23, 12, 15, 1; 17, 6, 3, 14, 20; 11, 22, 8, 19, 0.  Magic
\linebreak
$n \times n$ squares, such as this one, where the sums along the $n$-rows and columns, the
two diagonals and the $2n-2$ broken diagonals all have the same value are called
{\it pandiagonal} or ``diabolic" magic squares [{\bf 9}, ch.$~10$],[{\bf 11}].
To simplify, we shall use the term {\it ROW}
in all of the following to mean any one of a row, a column,
a diagonal or a broken diagonal.  Thus in a pandiagonal magic
square the sum along any ROW is the same and given by (1.1).

By a magic {\bf pandiagonal cube} of order $n$
we shall mean
an arrangement of the integers
$0, 1, 2,...n^3$ - 1 on the lattice points of an $n \times n \times n$ cube so that
the 3n squares parallel to the faces of the cube as well as the six
``diagonal" squares
which bisect the cube and contain its 4 body-diagonals are all pandiagonal magic squares.
Thus in each of the $3n$ + 6 (overlapping) squares contained in the cube, the sums along
all of the the ROWS---which here includes also the files---are
the same.
It is easily seen that in this case the common sum $\sigma_3$ must be
$$ \sigma_3 = \frac{n(n^3-1)}{2}
\eqno{(1.2)}
$$
\vspace{4pt}

By extension, we define a four dimensional (4-D) pandiagonal magic
\linebreak
\noindent
$n \times n \times n \times n$
cube as one in which the integers 0, 1, 2,...$n^4 - 1$ are placed, without
repeats, at the $n^4$ lattice sites of a 4-D hypercube, so that the three
dimensional (overlapping) $n \times n \times n$ cubes that can be formed within it,
are pandiagonal magic
cubes.
The sum along each ROW of a 4-D cube is easily shown to be
$$
\sigma_4 = n(n^4-1)/2
\eqno{(1.3)}
$$
%The important problem, that is the main subject of this paper, involves
%the algorithm to be used in constructing these cubes.
\vspace{4pt}

In constructing pandiagonal magic squares, cubes, and hypercubes, we follow the
idea of Euler [{\bf 6}]---who was concerned exclusively with
magic {\it squares}---and
used {\bf latin squares} (LS) in their construction.
As defined by  Denes and  Keedwell [{\bf 4}] and by Laywine and
Mullen [{\bf 9}]
a latin square of order $n$ is an $n \times n$ array each of whose
lattice points
is occupied by one of $n$ given symbols, in a way so that no row or column
contains any one of these symbols more than once.
The number of LS's grows rapidly
with $n$ and for example for $n = 10$ and 15
this number is $\sim 10^{36}$ and
\linebreak
$\sim 1.5 \times 10^{86}$, respectively [{\bf 9}, p.$~5$]. Here we are
interested in the relatively small subset of the LS's,
the {\bf pandiagonal latin squares},
for which in addition each diagonal as
well as each of the $2n-2$ broken diagonals---in
short it's ROWS---also contains each of the $n$ elements
precisely once.
%We call these {\bf pandiagonal latin squares}.
\vspace{4pt}

For a general algebraic theory of pandiagonal (diabolic) magic squares
see the analysis by Rosser and Walker [{\bf 11}].
The books by Andrews [{\bf 2}], Kraitchik [{\bf 5}] and
Rause Ball and Coxeter [{\bf 3}] contain
a more empirical approach to magic squares and include
various
examples and practical rules for constructing magic squares.
Pasles [{\bf 10}] describes
some unusual magic squares constructed
by Benjamin Franklin. %see the article by Pasles [10].

%Magic squares and cubes have been studied for many years and their status
%early in the 20th Century is summarized in the book by W. S. Andrews.[2]
%Martin Gardner [7]
%describes
Some 5 $\times$ 5
pandiagonal magic squares and some of their important features
are described by Gardner [{\bf 7}].  In a companion
column [{\bf 8}] he details a magic cube of order seven and notes its {\it non-pandiagonal
nature}.  This is
reenforced by the results of Wynne [{\bf 12}] who studied magic cubes of order seven
and showed that even if every square of the cube that is
parallel to a cube-face is pandiagonal,
not all of the six, diagonal squares of the cube can be pandiagonal.
This is consistent with the present
analysis, according to
which  only for $n \geq 11$ do magic pandiagonal cubes
exist.
(Similarly we find that only for $n \geq 17$ can 4-D pandiagonal
magic cubes exist.)
That {\bf non}-pandiagonal magic cubes of order 7 do exist, however, is
established by Andrews [{\bf 2}] and by Alspach and Heinrich [{\bf 1}],
the latter, incidental
to their discussion of cubes of order 4m.
\vspace{4pt}

%Much of the literature on magic squares and cubes of odd order is of a
%heuristic nature and not easily
%generalized to large values of n nor to hypercubes.  By contrast to these, the
%present approach makes use of an algorithm which, in principle, is able to encompass
%arbitrarily high odd values of n and can
%in principle
%be extended to arbitrarily high dimensions.
%\vspace{5pt}

\noindent
{\bf II. Odd-$n$ Pandiagonal Latin Squares and Cubes in 3 and 4 Dimensions
}

\vspace{15pt}

\noindent
{\bf A. Pandiagonal Latin Squares}
\vspace{8pt}

To set the stage for our discussion of magic pandiagonal cubes
in 3 and 4 dimensions, in this section we collect---and in some
instances amplify on---properties of pandiagonal latin squares [{\bf 4}],[{\bf 9}].

A {\bf latin square} of odd order $n$, {\bf LS},
is an $n\times n$ array involving $n$
distinct symbols with each symbol appearing once and
only once in each row and column.
We shall invariably take the $n$ symbols to be the
integers comprising the set $S$ defined by
$$  S = \{0,1,2,...,n-1\}
   \eqno{(2.1)}
$$
Also useful will be the set $\bar{S}$ defined by
$$  \bar{S} = \{1,2,...,n-1\}
   \eqno{(2.2)}
$$
A {\bf diagonal LS} has the additional property that each
of its two diagonals also contains the $n$ chosen symbols exactly once.
Finally a {\bf pandiagonal LS} is a diagonal LS in which, in addition,
each of the
$n$ symbols appears once and only once in each of the $2n-2$
``broken
diagonals", as defined in the preceding section.

Figure II shows 4 LS's of order 3,4,5, and 5, respectively.
The $~3 \times 3~$ LS in
part (a) is not a diagonal LS by virtue of
the three zeros

$$ \everymath={\displaystyle}
\begin{array}{c}
~~~~~~~~
\begin{array}{ccc}
       0 & 1 & 2 \\
       2 & 0 & 1 \\
       1 & 2 & 0 \\
         & (a) & \\
    \end{array}
         ~~~~~~~~~~~~~~~~~~
    \begin{array}{cccc}
           0 & 1 & 2 & 3 \\
           2 & 3 & 0 & 1 \\
           3 & 2 & 1 & 0 \\
           1 & 0 & 3 & 2 \\
             & (b) & & \\
         \end{array}
            ~~~~~~~~~~~
\\
    \begin{array}{ccccc}
        3 & 4 & 0 & 1 & 2\\
        1 & 2 & 3 & 4 & 0\\
        4 & 0 & 1 & 2 & 3 \\
        2 & 3 & 4 & 0 & 1 \\
        0 & 1 & 2 & 3 & 4 \\
          &   & (c) & & \\
         \end{array}
       ~~~~~~~~~~

     \begin{array}{ccccc}
      4 & 0 & 1 & 2 & 3\\
      2 & 3 & 4 & 0 & 1 \\
      0 & 1 & 2 & 3 & 4 \\
      3 & 4 & 0 & 1 & 2 \\
      1 & 2 & 3 & 4 & 0 \\
        &   & (d) & & \\
      \end{array}
  \end{array}
$$

\begin{center}
{\bf Figure II}
\end{center}

$~$
\vspace{10pt}

\noindent
appearing in the upper
left to lower right diagonal.
By contrast the $4\times 4$ LS in part (b) is diagonal but not pandiagonal
while the $5\times 5$, LS in part (c) is pandiagonal.
Note
that, given a pandiagonal LS, if we permute its symbols
$\{0,1,2,...n-1\}$ it remains pandiagonal.
For example, if in IIc, we carry out the permutation
0 $\rightarrow$ 1;
1 $\rightarrow$ 2;
2 $\rightarrow$ 3,
3 $\rightarrow$ 4;
\linebreak
4 $\rightarrow$ 0,
the resultant LS remains pandiagonal and is given by Figure IId.
Thus a given $n \times n$ pandiagonal LS is the basis for
$n$! different ones that result from the $n!$ possible permutations
of the $n$-symbols.

Consider for odd $n$, the $n\times n$ array $L$ with elements $L_{ij}$

$$  L_{ij} \equiv \alpha_1 i + \alpha_{2}j~~~~(mod~n)
   \eqno{(2.3)}
$$
where $i,j$ run over the elements of $S$
as do the elements $L_{ij}$ themselves
and $\alpha_1$ and $\alpha_2$
are non-zero positive integer parameters and thus
are elements of $\bar{S}$ in (2.2).
When we wish to stress the dependence of $L$ on
$\alpha_1$ and $\alpha_2$,
we shall also use the notation
$L =~ <\alpha_1, \alpha_2>$.

We now establish:

\noindent
{\bf Theorem (2.1)}:  If $L =~ <\alpha_1,\alpha_2>$
is the array in (2.3) and the greatest common divisor with $n$
of each of $\alpha_1, \alpha_2, \alpha_1 \pm \alpha_2$ is
1, that is
$$  (\alpha_i,n) = 1; ~~~~~i = 1,2  \eqno{(2.4a)}
$$
$$  (\alpha_1 +\alpha_2, n) = 1 \eqno{(2.4b)}
$$
$$  (\alpha_1 - \alpha_2, n) =1 \eqno{(2.4c)}
$$
then $L$ is a pandiagonal LS for $n \geq 5$.
(Note that (2.4) cannot be satisfied for even $n$.)

\noindent
{\bf Proof:}  Firstly, since for $n=3$, the possible $\alpha$-values
are 1 and 2 and since these violate (2.4b), it follows that (2.3) is not a
pandiagonal LS for $n=3$.  On the other hand,
for $n=5$, the pairs $[\alpha_1,\alpha_2] = [1,2]$ and
$[1,3]$, for example, do satisfy each of (2.4).

Secondly to establish first that $L$ in 2.3 for $n \geq 5$ is a diagonal LS, consider it
for fixed $j$, say, as $i$ ranges  over the
$j$th row.  As $i$ thus
varies from 0
\linebreak
to $n-1$, the $j$th row of $L$ varies over
the same set by virtue of the hypothesis $(\alpha_1,n) = 1$.
Similarly for fixed $i$, as $j$ varies over the $i$th column,
since $(\alpha_2, n) = 1$,
$L_{ij}$ varies over the same set S. %{\bf without repeats}.
Thus the rows and columns of $L$ satisfy the condition that $L$ be a
latin square.  Along the diagonal
from $(i,j) = (0, 0)$ to $(n - 1, n - 1)$,
$i=j$, so that along here
$L_{ij} ~\equiv ~(\alpha_1~ + ~\alpha_2)i$
($mod$ $n$).
Thus
again as $i$ varies over the set S the $n$ diagonal elements of $L$
must be some permutation
of S since $(\alpha_1 + \alpha_2,n) = 1$ according to (2.4b).  Finally
along the other
the diagonal from $(i,j) = (0, n - 1)$ to $(n - 1, 0)$, $i+j = n-1$,
so that $L_{ij} \equiv (\alpha_1 - \alpha_2)
i + \alpha_2(n-1)$ mod $n$ which has the same property by virtue of
(2.4c), $(\alpha_1 - \alpha_2, n) =1 $.  Thus we have established that
$L$ in (2.3) under the constraints in (2.4)
is a diagonal LS; it remains only to establish that it is also
pandiagonal.

To this end, consider, for example, the ``split"
diagonal just above---and parallel to---the lower
left to upper right diagonal of $L$
and its appendage in the lower right hand corner,
$(i,j) = (n-1,0)$.
(We assume $i$ increases from 0 to the right and $j$ increases
upward from 0.)
Its entries are defined by
$j=i+1$,
$i=0,1,2,...n-2$ and
$(i,j) = (n-1, 0)$.
Substituting $j= i+1$ into (2.3) we find
$$  L_{i, i+1} \equiv (\alpha_1 + \alpha_2)i + \alpha_2~~ (mod~n)
$$
If we now let $i$ run over the complete set S,
it is easy to see that
since $L_{in} \equiv L_{io}$ $(mod~n)$
we obtain
both parts of the
broken diagonal and that no two of these $n$-elements are the same
since they are simply a permutation of the elements of S
by virtue of (2.4b).
This argument is easily repeated for all broken diagonals and we
conclude that each broken diagonal consists of a
permutation of the elements of  S.
We have thus established the theorem{\bf .}
%that L as defined in (2.3), subjected to the constraints
%in (2.4), is a pandiagonal latin square.
%%=====================

\noindent
{\bf Remark 2.1:}  If $\alpha_1,\alpha_2$ satisfy each
of (2.4), then so do $k\alpha_1, k\alpha_2$ for any $k \epsilon \bar{S}$
for which $(k,n) = 1$.

\noindent
{\bf Remark 2.2:}  Reference to (2.3) shows that if
$L = <\alpha_1,\alpha_2>$ is a pandiagonal LS, then
so is $<k\alpha_1, k\alpha_2>$ for any positive integer $k\epsilon\bar{S}$
provided $(k,n) = 1$.  For according to (2.3),
$<k\alpha_1, k\alpha_2>$, is simply
$<\alpha_1,\alpha_2>$ with its elements permuted in
some way.
%Similarly, $\{\alpha_1,\alpha_2\} + x~(mod~n)$ where $x$
%is an integer not divisible by $n$, is simply
%$\{\alpha_1,\alpha_2\}$ with its elements permuted in some way.

\noindent
{\bf Remark 2.3:}  If we add an integer $x \epsilon \bar{S}$ to the
right hand side of (2.3), we obtain $<\alpha_1, \alpha_2>$ with its
elements permuted in some way.

An important notion relating to LS's is that of the
{\bf orthogonality} [{\bf 9}, ch. 2].
Two
$n\times n$,
LS's are said to be orthogonal, if when
they are superposed, none of the
$n^2$
ordered pairs of
elements that result occurs more than once.  Thus if
$L_{ij}^{(1)}$
and
$L_{ij}^{(2)}$
$\epsilon S$
are any corresponding elements of the two LS's,
$L^{(1)}$ and
$L^{(2)}$
then they are orthogonal if and only
if the ordered pairs
$[L_{ij}^{(1)},
L_{ij}^{(2)}]$  and
$[L_{\alpha\beta}^{(1)},
L_{\alpha\beta}^{(2)}]$
differ
$(mod~n)$ for any fixed values of
$i,j \epsilon S$
but for all choices of
$\alpha, \beta, \epsilon S$.

\noindent
{\bf Remark 2.4}  If two LS's are orthogonal,
they remain so when the elements
of either or both undergo arbitrary permutations.

%Consider two pandiagonal LS's, L$^{(1)}  = \{\alpha_1,
%\alpha_2\}$ and $L^{(2)} = \{\beta_1,\beta_2\}$.
We now establish:

\noindent
{\bf Theorem (2.2)}  Let $L^{(1)} = <\alpha_1,\alpha_2>$
and $L^{(2)} =~<\beta_1,\beta_2>$ be two pandiagonal
LS's with both pairs $\alpha_1,\alpha_2$ and $\beta_1,\beta_2$
each satisfying (2.4).
If the determinant $d_2$ defined by

$$  d_2 = \Biggr| \begin{array}{cc}    \alpha_1 & \alpha_2\\
                                           \beta_1 & \beta_2 \end{array} \Biggr|
   \eqno{(2.5)}
$$
is relatively prime to $n$,
i.e.

$$   (d_2,n)=1
  \eqno{(2.6)}
$$

\noindent
then $L^{(1)}$ and $L^{(2)}$ are orthogonal
pandiagonal
LS's.

\noindent
{\bf Proof:}  Suppose for $i,j \epsilon S$ there is a second
pair $k,\ell \epsilon {S}$ for which the ordered pairs
$[L_{i}j^{(1)}, L_{i}j^{(2)}] \equiv [L_{k\ell}^{(1)}, L_{k\ell}^{(2)}]$
mod $(n)$.
That is, suppose

\singlespace
$$  \begin{array}{rcll}
     \alpha_1 i + \alpha_{2}j & \equiv & \alpha_1 k+ \alpha_{2}\ell&\\
               & & & (mod~n)\\
     \beta_1 i  + \beta_{2}j & \equiv & \beta_1 k + \beta_2 \ell &
    \end{array}
%   \eqno{(2.10)}
$$

\doublespace

\noindent
which for convenience we express in matrix notation
$$  \Biggr( \begin{array}{cc} \alpha_1 & \alpha_2 \\
                              \beta_1  & \beta_2 \end{array} \Biggr)
     \left(        \begin{array}{l} (i-k) \\
                                       (j-\ell) \end{array} \right) \equiv 0
    ~~~~~(mod~n)
    %\eqno{(2.11)}
$$
Now since by hypothesis $(d_2,n) = 1$, so that
in particular
$d_2 \not\equiv 0$ $(mod~n)$, we can
multiply this relation by
$$  d_2 \Biggr( \begin{array}{cc}  \alpha_1 & \alpha_2 \\
                                        \beta_1  & \beta_2 \end{array} \Biggr)^{-1} =
             \Biggr( \begin{array}{rr} \beta_2 & -\alpha_2 \\
                                       -\beta_1 & \alpha_1 \end{array} \Biggr)
$$
so that it becomes
$$  d_2 \biggr( \begin{array}{c} (i-k) \\
                              (j-\ell) \end{array} \biggr) \equiv 0~~~mod~n
$$

\noindent
Finally, since $(d_2,n) = 1$ it follows
that
$i=k$ and
$j=\ell$.
Thus
$L^{(1)}$ and
$L^{(2)}$ are orthogonal pandiagonal LS's{\bf .}

The orthogonality criterion of the two LS's in (2.6) is very convenient
and as we shall see is extendable to higher dimensions.

Let us consider the question as to the number
of distinct pairs $[\alpha_1,\alpha_2]$ there are for given $n$
for which $L = <\alpha_1,\alpha_2>$ is a pandiagonal LS.  To
simplify let us assume in the following that $n$ is a prime
$p$.
For a given prime $p \geq 5$, consider the $p-3$, $[\alpha_1,\alpha_2]$,
pairs
$$  [\alpha_1,\alpha_2] = [1,2], [1,3], [1,4],... [1,p-2]
  \eqno{(2.7)}
$$
Obviously each of these satisfies (2.4) and no other pair for which
$\alpha_1 = 1$ does so.

%\pagebreak

\noindent
{\bf Theorem (2.3):}  The pandiagonal latin squares associated
with the $[\alpha_1,\alpha_2]$ pairs in (2.7) are:
\begin{description}
\item[(1)] pairwise mutually orthogonal
\item[(2)] any other pandiagonal LS can be obtained from one
associated with
$<1,\ell>$,  $\ell = 2,3,...p-2$,
in (2.7) by a permutation of symbols.
\end{description}

\noindent
{\bf Proof:}
\begin{description}
\item[(1)] Consider the two pandiagonal LS's $<1,\ell>, <1,m>;
\ell, m = 2,3,...p-2,$ $(\ell \not= m)$.
According to (2.5) the associated determinant $d_2$ is
$$  d_2 = \left| \begin{array}{ll} 1 & \ell \\ 1 & m .
\end{array} \right|
= m-\ell
$$
and obviously satisfies (2.6) since we assume $m\not= \ell$.

\item[(2)]  Let $x,y$ be any two unequal positive integers $\epsilon \bar{S} =
\{1,2,...p-1\}$ that satisfy each of (2.4) so that
$<x,y>$ is a pandiagonal LS.  If $x^{-1}$ is the inverse of
$x~(mod~p)$ then making use of (2.3) and Remark (2.2),
we find that
$<1,x^{-1}y>$ is also pandiagonal and is obtained from
$<x,y>$ by a permutation of its symbols.
Finally since $x^{-1}y~(mod~ p)$ must be one of
$(2,3,...p-2)$ in (2.7) it follows from Remarks 2.1 and
2.2 that $<x,y>$ can be obtained from $<1,x^{-1}y>$
by a permutation of its symbols{\bf .}
\end{description}

%===========================================

For example, for $p = 7$, successive multiplication $(mod~p)$ by use of
$k = 2,3,...6$,
leads to $<1,2>$ $\rightarrow$ $<2,4>, <3,6>, <4,1>, <5,3>,
<6,5>$. Similarly $<1,3>$
$\rightarrow$ $<2,6>$, $<3,2>$, $<4,5>$, $<5,1>$,
$<6,4>$; and $<1,4>$ $\rightarrow$ $<2,1>$,
$<3,5>$, $<4,2>$, $<5,6>$, $<6,3>$.

\noindent
{\bf Remark 2.5:}  If we allow for permutation of symbols
%in any given pandiagonal LS in (2.3)
then all pandiagonal
LS's in (2.3) can be obtained, for given $p$ by use only of
the $[\alpha_1,\alpha_2]$ pairs in (2.7).

Finally, as shown
for pandiagonal {\it magic squares}, by Ball and Coxeter [3, p. 203] and by
Martin Gardner [{\bf 7}] for a $5\times 5$ pandiagonal
{\it magic} square, we find that
pandiagonal LS's have an analogous unusual property.
If in an $n\times n$ pandiagonal LS we move the left hand
column so that it becomes the right hand column (or vice versa)
or similarly move the top row to the bottom, the resultant
array is again a pandiagonal LS.  This has the consequence
that if we ``tile" the plane with a given $n\times n$ pandiagonal LS we can
outline any $n\times n$ square on this infinite pattern and obtain
a pandiagonal LS.

%Adding paragraph 9/26/01:
The underlying result is contained in:

\noindent
{\bf Theorem 2.4:}  If the left column of a pandiagonal latin square
is moved so it becomes the right column, the resulting LS is a pandiagonal
LS obtained from the original one by a permutation of symbols.

\noindent
{\bf Proof:}  Let $L = <\alpha_1,\alpha_2>$ be the original pandiagonal LS
and define a second one $L^\prime$ with elements given by
$$  L_{ij}^\prime \equiv L_{ij} + \alpha_1 \equiv (\alpha_1 + 1)i +
     \alpha_2j~~~~(mod~n)
$$
It follows from (2.3) that $L^\prime$ is simply $L$ with its elements
permuted.  Further we have

%To show that $L = \{\alpha_1,\alpha_2\}$ in (2.3) with
%the constraints in (2.4) has this property, consider a new
%latin square $L^\prime \equiv \{\alpha_1,\alpha_2\} + \alpha_1~(mod~n)$.
%This square has the elements:
%$$  L_{oj}\equiv \alpha_1i + \alpha_2j + \alpha_1 \equiv \alpha_1(i+1) + \alpha_2j
%    ~~~~(mod~n)
%$$
%Obviously
$$L_{oj}^{\prime}\equiv L_{1j}; L_{1j}^\prime \equiv
   L_{2j};... L_{n-1,j}^\prime \equiv L_{oj}; (mod~n);~j=0,1,...n-1
$$
which shows that $L^\prime$ is simply $L$ with its left column
moved so it becomes the right column.
%which establishes the result.
%Also by Remark 2.2 we note that $L^\prime$ is simply $L$
%with its elements permuted in some way.
%Similarly it is easily shown
%that the top rows can be moved to the bottom one without
%changing the pandiagonal nature of the LS.

For example, if the left hand column of the LS in Figure IIc,
which incidentally is simply $<1,3>$, (with $i$ increasing to the right and
j increasing upwards) is moved to the right side, the original LS,
but with permuted symbols,
in IId results!
Of considerable interest perhaps is that a form of this property, as
will be shown below, has an analogue in higher dimensions.

%=======================================
%%end{here} 6/13/01

\noindent
{\bf B. Pandiagonal Latin Cubes}

We define a {\bf pandiagonal latin cube} as an $n\times n\times n$ cube
each of whose $n^3$ lattice points contains one of the members of
$S = \{0,1,2,...n-1\}$ and in a way so that each of its $3n+6$ constituent squares
is a pandiagonal latin square.   Recall in this connection that an
$n\times n\times n$ cube has $3n$ squares parallel to a cube face plus 6
``diagonal" squares which contain its 4 body diagonals.
Note
that, as here defined, in
a pandiagonal latin cube each row, column, file, diagonal
and broken diagonal of each of its squares, i.e., its ROWS, contain
each element of $S$ once and only once.

In an obvious generalization of (2.3) to three dimensions, we define an
$n \times n \times n$
array C by the formula
$$ C_{ijk} \equiv \alpha_1 i + \alpha_2 j + \alpha_3 k~~~~(mod~n)
\eqno{(2.8)}
$$

\noindent
where each element of $C_{ijk}$ $\epsilon~S, \{0, 1,...n - 1\}$.
Here {\it i,j,k} are integer variables each running over
$\{0,1,...n-1\}$, and $\alpha_1,
\alpha_2, \alpha_3$ are elements
of
\linebreak
$\bar{S} = \{1, 2,...n - 1\}$.
We shall use the notation $C =~ <\alpha_1,\alpha_2,\alpha_3>$
when we wish to focus on the dependence of $C$ on
$\alpha_1,\alpha_2,\alpha_3$.
We now establish the following:
\vspace{8pt}

\noindent
{\bf Theorem (2.5):} The cube defined in (2.8) is a pandiagonal
latin cube---in that each of its constituent 3n + 6
squares is pandiagonal---provided $\alpha_1, \alpha_2, \alpha_3$
satisfy the constraints:

\singlespace

$$ (\alpha_\ell, n) =  1 ;~~~~ \ell = 1, 2, 3
\eqno{(2.9a)}
$$

$$ (\alpha_\ell + \alpha_{\ell^\prime}, n) = 1;~~~~ \ell, \ell^\prime = 1, 2, 3;~ \ell \not= \ell^{\prime}
\eqno{(2.9b)}
$$

$$ (\alpha_\ell - \alpha_{\ell^\prime}, n) = 1;
   ~~~~ \ell, \ell^\prime = 1, 2, 3;~  \ell \not= \ell^\prime
\eqno{(2.9c)}
$$

$$ (A, n) = 1 ; ~~~~ A = \alpha_1 + \alpha_2 + \alpha_3
\eqno{(2.9d)}
$$

$$ (A - 2 \alpha_\ell, n) = 1; ~~~~ \ell = 1, 2, 3
\eqno{(2.9e)}
$$

\doublespace

\noindent
{\bf Proof:} Consider first the
$3n$ squares parallel to the faces of the
cube. For fixed $k$, say, 0 $\leq k \leq n - 1$,
consider the square of C$_{ijk}$ that is parallel to the
$i-j$ plane as $i$ and $j$ vary over
the elements of S.  Then C$_{ijk}$ has essentially
the same structure as does $L_{ij}$ in
(2.3) since the constant $\alpha_3k$ is of no consequence.
Making use of the restrictions on the
$\alpha$'s in (2.9 a,b,c) and comparing with
those in (2.4 a, b, c)
we conclude that for any fixed $k$, the square $C_{ijk}$  is a
pandiagonal latin square.  Repeating this argument
for fixed $i$, with $j$ and $k$ variable and for fixed
$j$ with $i$ and $k$ variable we conclude that all 3n squares in the cube defined by
(2.8) that are parallel to a cube face are pandiagonal latin squares.
\vspace{8pt}

%Consider now the six squares not parallel to a cube face.  For the ``diagonal
%square" with upper edge D$_1$ in Figure V, i = j so that we have for this square
%$$ C_{iik} \equiv (\alpha_1 + \alpha_2) i + \alpha_3k~~(mod~n)
%$$

With regard to the six squares not parallel to a cube
face we proceed as follows.  For that diagonal square, with
vertices at $(i,j,k) = (0,0,0),
\linebreak
(0,0,n-1),(n-1,n-1,n-1),(n-1,n-1,0)$,
we have $i=j$ so that (2.8) becomes
$$  C_{iik} \equiv (\alpha_1 + \alpha_2) i+\alpha_3k~~~~~~~~(mod~n)~~.
$$
But this is again of the form in (2.3) if we make the replacements
in (2.4a,b,c) $\alpha_1 \rightarrow \alpha_1 + \alpha_2;
\alpha_2 \rightarrow \alpha_3$.  The first of (2.4) is satisfied
because of (2.9a,b) while (2.4b) and (2.4c) become respectively
$(\alpha_1 + \alpha_2 + \alpha_3,n) =1$ and
$(\alpha_1 + \alpha_2-\alpha_3,n)=1$
and these are the same constraints as in (2.9d) and (2.9e),
respectively.  Thus $C_{iik}$, the given square, is a
pandiagonal latin square.

Similarly for the diagonal square perpendicular to
$C_{iik}$ whose vertices have the coordinates
$(i,j,k) = (n-1,0,0), (0,n-1,0), (0,n-1,n-1),
\linebreak
(n-1,0,n-1)$, we have $i+j = n-1$  and

$$  C_{i,(n-1-i),k} \equiv (\alpha_1-\alpha_2) i
     + \alpha_3k + \alpha_2(n-1)~~~~~~(mod~n)
$$

\noindent
which on comparison with (2.3) and (2,4) with the replacements
$\alpha_1 \rightarrow
\alpha_1 - \alpha_2; \alpha_2 \rightarrow \alpha_3$
leads to the conditions
$ (\alpha_1 - \alpha_2 + \alpha_3, n ) = 1
$,
and
$ (\alpha_1 - \alpha_2 - \alpha_3, n) = 1
$
The first of these is the same as (2.9e) with $\ell$ = 2
and the second is the
same as (2.9e) with $\ell$ = 1 if we make use of the
fact that (-x, y) = 1 is equivalent
to (x,y) = 1.

A similar argument shows that for the constraints in (2.9)
the remaining four ``diagonal
squares" of the cube are also pandiagonal latin squares.
The theorem is thus established{\bf .}

Just as for the 2-D case, we require a 3-D analogue of orthogonality
of latin cubes.  For our purposes, we shall say that three latin
cubes are {\bf orthogonal} if when
they are superposed none of the $n^3$ ordered triplets of
elements that result occurs more than once.  There are
other definitions of orthogonal latin cubes, [{\bf 9}, ch. 3],
but for purposes of producing
{\it magic} cubes this definition is essential. %the most useful.
\vspace{10pt}

We now establish:

\noindent
{\bf Theorem (2.6):}  Consider the 3 pandiagonal latin cubes
$$  C_{ijk}^{(q)} \equiv \alpha_{1q} i + \alpha_{2q}j
     + \alpha_{3q} k~~~~~~~~~(mod~n);~~~~~~~~~~q = 1,2,3
   \eqno{(2.10)}
$$
where $\alpha_{pq}$ $(p,q = 1,2,3)$ are elements of
$\bar{S}$ and for each value for $q$ satisfy the conditions
in (2.9) and let $d_3$ be the determinant
$d_3 = |\alpha_{pq}|$.  Then if
$d_3$ is relatively prime to $n$, that is
$$  (d_{3},n) = 1
  \eqno{(2.11)}
$$
then the three cubes $C^{(q)}$, $q=1,2,3$ are orthogonal.
\vspace{10pt}

\noindent
{\bf Proof:}  Suppose on the contrary that for a given $(i,j,k)$
there existed an integer triplet $(u,v,w)$ each $\epsilon S$
not equal to $(i,j,k)$ for which
$C_{ijk}^q \equiv C_{uvw}^q$ $(mod~n)$ for each $q$.  Then,
as for the 2-D case, we could express this in matrix
notation
$$   \left( \begin{array}{rrr}
    \alpha_{11} & \alpha_{21} & \alpha_{31}\\
    \alpha_{12} & \alpha_{22} & \alpha_{32} \\
    \alpha_{13} & \alpha_{23} & \alpha_{33} \end{array} \right)
   \left( \begin{array}{l}
       i-u \\ j-v\\ k-w \end{array} \right) \equiv 0~~~~~~mod~n
   \eqno{(2.12)}
$$
Since by (2.11), the determinant  $d_3$ is relatively prime to $n$, so that
in particular
$d_3 \not\equiv 0$ $mod~n$, we may multiply both sides of (2.12)
by $d_3$ times the inverse of the matrix on the left.
The result is
$$  d_3 \left( \begin{array}{l}
      i-u\\  j-v\\  k-w \end{array} \right) \equiv 0
~~~~~~mod~n.
$$
Finally since $(d_3,n) = 1$ it follows that $i=u; j=v; k=w$
and the theorem is proved{\bf .}

An empirical study of the constraints in (2.9) shows that no triplets
$[\alpha_1,\alpha_2, \alpha_3]$ satisfying (2.9) exists for
$n=5,7,9$.  (The latter clearly since
$n$ would be divisible
by 3.)   Such an analysis is most easily carried out by
recognizing
that---as in the 2-D case---without loss of generality we
can take
$\alpha_1 = 1$ and consider simply the triplets
$[1,\alpha_2,\alpha_3]$.  For
$n=7$, it is easily confirmed that
no values, such
$[1,2,3], [1,2,4], [1,2,5]$ would satisfy all
of (2.9).  For
$n=11$ however we find, among others, the
possibilities
$[1,2,4], [1,2,5], [1,2,6], [1,2,7], [1,5,8]$,
$[1,6,8]$
as well as these with
$\alpha_{2}$ and
$\alpha_3$ interchanged.
As an example of {\it orthogonal} pandiagonal cubes we
note that for appropriate integers,
$\ell, m, p $
for the three LS's
$<1,2,\ell>$,
$<1,m,2>$ and
$<1,2,p>$,
for
$n=11$,
$d_3$ is given by
\vspace{10pt}

$$  d_3 = \left| \begin{array}{lll}
      1 & 2 & \ell \\
      1 & m & 2 \\
      1 & 2 & p \end{array} \right| =
       (m-2)(p-\ell),
$$
\vspace{10pt}

\noindent
and will for values of $\ell,m,p$ with
$\ell \not= p$ and
$4 \leq \ell,m,p \leq 9$
lead to $(d_3, n) = 1$.    Figure III shows two planar
sections through the cube $<1,2,7>$ for $n=11:$ (a)
corresponds to the square $k=2$ in (2.8) and (b) to
the diagonal square $i=j$.

Note, however, that for the triplet
$<1,2,\ell>, <1,2,m>, <1,2,p>$, $d_3 = 0$ so that
these three do not constitute orthogonal cubes as we have
defined them even though each cube is itself pandiagonal. %latin cube.

As for the analogous 2-D case (Theorem 2.4) we can easily establish the fact
that if we move, say,$~$ a face of a pandiagonal$~$ latin cube to its
opposite
\pagebreak

\kimspace

$$
\begin{array}{l}
\begin{array}{l}
\left[ \begin{array}{rrrrrrrrrrr}
3 & 5 & 7 & 9 & 0 & 2 & 4 & 6 & 8 & 10 & 1\\
4 & 6 & 8 & 10 & 1 & 3 & 5 & 7 & 9 & 0 & 2\\
5 & 7 & 9 & 0 & 2 & 4 & 6 & 8 & 10 & 1 & 3\\
6 & 8 & 10 & 1 & 3 & 5 & 7 & 9 & 0 & 2 & 4\\
7 & 9 & 0 & 2 & 4 & 6 & 8 & 10 & 1 & 3 & 5\\
8 & 10 & 1 & 3 & 5 & 7 & 9 & 0 & 2 & 4 & 6\\
9 & 0 & 2 & 4 & 6 & 8 & 10 & 1 & 3 & 5 & 7\\
10 &1 & 3 & 5 & 7 & 9 & 0 & 2 & 4 & 6 &  8 \\
0 & 2 & 4 & 6 & 8 & 10 & 1 & 3 & 5 & 7 & 9\\
1 & 3 & 5 & 7 & 9 & 0 & 2 & 4 & 6 & 8 & 10\\
2 & 4 & 6 & 8 & 10 & 1 & 3 & 5 & 7 & 9 & 0
\end{array} \right] \\
~~~~~~~~~~~~~~~~~~~~~~~~~~~~~~~(a)
\end{array}
\\
\begin{array}{l}
\left[ \begin{array}{rrrrrrrrrrr}
0 & 7 & 3 & 10 & 6 & 2 & 9 & 5 & 1 & 8 & 4\\
3 & 10 & 6 & 2 & 9 & 5 & 1 & 8 & 4 & 0 & 7\\
6 & 2 & 9 & 5 & 1 & 8 & 4 & 0 & 7 & 3 & 10\\
9 & 5 & 1 & 8 & 4 & 0 & 7 & 3 & 10 & 6 & 2\\
1 & 8 & 4 & 0 & 7 & 3 & 10 & 6 & 2 & 9 & 5\\
4 & 0 & 7 & 3 & 10 & 6 & 2 & 9 & 5 & 1 & 8\\
7 & 3 & 10 & 6 & 2 & 9 & 5 & 1 & 8 & 4 & 0\\
10 & 6 & 2 & 9 & 5 & 1 & 8 & 4 & 0 & 7 & 3\\
2 & 9 & 5 & 1 & 8 & 4 & 0 & 7 & 3 & 10 & 6\\
5 & 1 & 8 & 4 & 0 & 7 & 3 & 10 & 6 & 2 & 9\\
8 & 4 & 0 & 7 & 3 & 10 & 6 & 2 & 9 & 5 & 1
\end{array}\right] \\
~~~~~~~~~~~~~~~~~~~~~~~~~~~~~~~~(b)
\end{array}
\end{array}
$$

\hspace{2in}{\bf Figure III}
\vspace{10pt}

%\pagebreak

\doublespace

\noindent
side, the resultant cube remains a
pandiagonal latin cube.
To see this, consider
$C =~ <\alpha_1,\alpha_2,\alpha_3>$ in (2.8) and define
%$C^\prime =~ <\alpha_1,\alpha_2, \alpha_3> +
%    \alpha_1$
$C^\prime$ with elements given by $C_{ijk}^\prime = C_{ijk} + \alpha_1$
so that

$$  C_{ijk}^\prime \equiv \alpha_1(i+1) + \alpha_2 j + \alpha_3 k~(mod~n) .
$$

\noindent
Obviously,
$$  C_{ojk}^\prime = C_{1jk}; C_{1jk}^\prime \equiv C_{2jk};...
    C^\prime_{n-1,jk} \equiv C_{o,jk}~~~~~~~ (mod~n);
     ~~0 \leq j,k \leq n-1
$$
Thus $C^\prime$ is obtained from $C$ by transporting its
$i=0$ face to $i=n-1$.  This implies that if we ``tile"
all of 3-D space with a given pandiagonal latin cube {\bf any}
$n\times n\times n$ cube selected out of this infinite
array will be a pandiagonal latin cube
%Such a cube will be the original one
but with its elements permuted
in some way.

\noindent
{\bf C.  Four dimensional pandiagonal latin cubes}

A 4-D pandiagonal latin cube is an arrangement of the
integers
\linebreak
$\{0,1,2,...n-1\}$ among the cells of an
$n\times n \times n \times n$ cube in a way so that
each of its $4n+12$ constituent 3-D cubes is a
pandiagonal latin cube.

In an obvious generalization of (2.8), we define
an $n \times n \times n \times n$ array of integers by the formula
$$ H_{ijk \ell} \equiv \alpha_1i + \alpha_2j + \alpha_3k +
    \alpha_4 \ell~~~~~~~~(mod~n)
\eqno{(2.13)}
$$

\noindent
where each of $H_{ijk\ell}$ is $\epsilon$ $S$ = $\{0,1,2,...,n-1\}$
and where $i,j,k$, $\ell$
are integer variables each with the range $0,1,2,...,n-1$.  The four
quantities
$\alpha_1, \alpha_2, \alpha_3, \alpha_4$ are integer parameters
$\epsilon~\bar{S} = \{1,2,...n-1\}$.
As before we use the notation $<\alpha_1,\alpha_2,\alpha_3,\alpha_4>$
when the dependence of $H$ on the $\alpha$'s is of interest.
We now establish the following:
\vspace{8pt}

\noindent
{\bf Theorem (2.7):} The hypercube defined in (2.13)
is a 4-D pandiagonal latin cube provided the four
integer parameters $\alpha_1, \alpha_2, \alpha_3, \alpha_4$ satisfy the constraints
\kimspace
$$
\begin{array}{llr}
(\alpha_m, n) = 1 & ;m = 1,2,3,4 & (2.14a) \\
(\alpha_m + \alpha_{m^\prime}, n) = 1 & ;m, m^\prime = 1,2,3,4; m \not= m^\prime & (2.14b) \\
(\alpha_m - \alpha_{m^\prime}, n) = 1 & ;m, m^\prime = 1,2,3,4; m \not= m^\prime & (2.14c) \\
(B - \alpha_{m}, n) = 1 & ;B = \alpha_1 + \alpha_2 + \alpha_3 + \alpha_4; m = 1,2,3,4 &
(2.14d) \\
(B,n) = 1 & ; & (2.14e) \\
(B - 2\alpha_m, n) = 1 & ;m = 1, 2, 3, 4 & (2.14f) \\
(B - \alpha_{m^\prime} - 2\alpha_m,n) = 1 & ;m, m^\prime = 1,2,3,4; m \not= m^\prime & (2.14g) \\
(B - 2\alpha_{m^\prime} - 2\alpha_m,n) = 1 & ;m,m^\prime = 1,2,3,4; m \not= m^\prime & (2.14h)
\end{array}
$$

\doublespace

\noindent
{\bf Proof:} Consider first the 3-D cube that results for
fixed $\ell$ from (2.13) as $i,j,k$
vary over S.  It will be a pandiagonal latin cube provided the
constraints in (2.9a)-(2.9e)
are satisfied.  Now (2.9a)-(2.9e)
are contained within (by
appropriate choice of the subscripts)
(2.14a), (2.14b), (2.14c), (2.14e) and (2.14g),
respectively.
Thus since the added constant $\alpha_4\ell$ plays no role
the cube is a pandiagonal latin cube.
Similarly, the cubes that result
for fixed $i$, as $j,k, \ell$ vary, and for
fixed $j$ as $i,k,\ell$ vary and for fixed $k$ as $i,j, \ell$ vary over S
are all pandiagonal latin
cubes.  There are altogether
$4n$ pandiagonal 3-D latin cubes of this type.
\vspace{8pt}

Similarly the ``diagonal" 3-D cube that results for $i = j$, as $i, k, \ell$
vary over S has the
form
$$ H_{iik \ell} \equiv (\alpha_1 + \alpha_2)i +
    \alpha_3k + \alpha_4 \ell~~(mod~n)
$$
On comparison with (2.8) and (2.9a) - (2.9e) we
see that this is also a pandiagonal latin cube
since we assumed
(2.14a), (2.14b) (2.14c), (2.14d), (2.14e), (2.14g) to be satisfied .
And similarly for
the other five pairs: $i = k, i = \ell$, $j = k, j = \ell$, $k = \ell$.
\vspace{8pt}

Finally for the ``diagonal" latin 3-D cube that results
for $i + j = n - 1$, as $i, k$ and $\ell$
vary over S
$$H_{i, n-1-i,k,\ell} \equiv (\alpha_1 - \alpha_2)i + \alpha_3k
   + \alpha_4 \ell
   - \alpha_2 (n-1)~(mod~n)
$$

\noindent
which on comparison with (2.9a) - (2.9c) is a
pandiagonal latin cube by virtue of (2.14a), (2.14b),
(2.14c), (2.14f), (2.14g), (2.14h).
And similarly for the remaining five cubes
$i + k = n - 1, i + \ell = n - 1, j + k = n - 1,
j + \ell = n - 1$ and $k + \ell = n - 1$.
This proves the theorem{\bf .}
\vspace{8pt}

Turning to the question of orthogonality, we define four,
4-D pandiagonal cubes to be {\bf orthogonal}, if when they are
superposed no two of the $n^4$ ordered quartets of elements that result
are the same.  As in the lower dimensional cases there are other
definitions of this orthogonality [{\bf 9}], but for our purposes this
one is essential.

By analogy to Theorems (2.2) and (2.6) we have:

\noindent
{\bf Theorem (2.8):}  Consider the four, 4-D pandiagonal latin
cubes
$$   H_{ijk\ell}^{(q)} \equiv \alpha_{1q} i + \alpha_{2q}j
    + d_{3q} k + \alpha_{4q} \ell~~~~~(mod~n);~~~~~~
     q = 1,2,3,4
   \eqno{(2.15)}
$$
where $\alpha_{pq}~(p,q = 1,2,3,4)$ are elements of
$\bar{S}$ which for each value of $q$ satisfy
the conditions in (2.14) and the variables $i,j,k,\ell$
range over
\linebreak
$S = \{0,1,2,...n-1\}$.  Further, let
$d_4 = |\alpha_{pq}| $ be the determinant of the
$\alpha$'s.  Then if $d_4$ is relatively prime to $n$
$$ (d_4,n) = 1
   \eqno{(2.16)}
$$
then the four 4-D cubes in (2.15) are orthogonal.

\noindent
{\bf Proof:}  Suppose to the contrary there existed a quartet
of integers $(u,v,w,x)$ each $\epsilon$ $S$, not equal to
any given quartet $(i,j,k,\ell)$ for which $H_{ijk\ell}^{(q)}
\equiv H_{uvwx}^{(q)}$ $(mod~n)$ for each $q$.
Then just
as in deriving (2.12) we would find
$$  \left(  \begin{array}{rrrr}
     \alpha_{11} & \alpha_{21} & \alpha_{31} & \alpha_{41}\\
     \alpha_{12} & \alpha_{22} & \alpha_{32} & \alpha_{42}\\
     \alpha_{13} & \alpha_{23} & \alpha_{33} & \alpha_{43}\\
     \alpha_{14} & \alpha_{24} & \alpha_{34} & \alpha_{44}
    \end{array} \right)
    \left( \begin{array}{l}
       i-u \\ j-v\\ k-w\\ \ell - x \end{array} \right) \equiv 0
      ~~~~~~~~~~mod~n
    \eqno{(2.17)}
$$
and conclude following essentially the same steps as before that
$i,j,k,\ell$ must be equal to $u,v,w,x$ respectively.
Thus concluding the proof of orthogonality{\bf .}

It is easily confirmed by enumerating the various
possibilities that only for n $\geq$ 17 is
it possible to find
integers
$[\alpha_1,\alpha_2,\alpha_3,\alpha_4]$
among $\{1,2,...n-1\}$
for which (2.14) can be satisfied.  In particular for $n = 17$,
possible 4-D
latin hypercubes are given by
$<1,2,4,8>$, $<1,2,4,9>$, $<1,2,13,8>$, \linebreak $<1,2,13,9>$ as is readily
confirmed.  Thus a possible form for
the determinant
of the matrix in (2.17) is
$$
d_4 =
\left|
\begin{array}{cllclclclcl}
& 1 & 2 & 4 & 8 \\
& 1 & 2 & 4 & 9 \\
& 1 & 2 & 8 & 4   \\
& 1 & 4 & 9 & 2 \\
\end{array}\right| = -8
$$

\noindent
and since (-8, 17) = 1,
the four 4-D pandiagonal latin cubes
$<1,2,4,8>, \linebreak <1,2,4,9>, <1,2,8,4>, <1,4,9,2>$ constitute
an orthogonal set of such hypercubes, for $n=17$.
\vspace{4pt}

\kimspace
Figure IV, shows a planar section through the hypercube
$<1,2,4,9>$ corresponding to
$i = 2, j + k = 16$ in (2.13).

$$
\begin{array}{l}
\left[ \begin{array}{rrrrrrrrrrrrrrrrr}
0  & 9  & 1  & 10 & 2  & 11 & 3  & 12 & 4  & 13 & 5  & 14 & 6  & 15 & 7  & 16 & 8  \\
2  & 11 & 3  & 12 & 4  & 13 & 5  & 14 & 6  & 15 & 7  & 16 & 8  & 0  & 9  & 1  & 10 \\
4  & 13 & 5  & 14 & 6  & 15 & 7  & 16 & 8  & 0  & 9  & 1  & 10 & 2  & 11 & 3  & 12 \\
6  & 15 & 7  & 16 & 8  & 0  & 9  & 1  & 10 & 2  & 11 & 3  & 12 & 4  & 13 & 5  & 14 \\
8  & 0  & 9  & 1  & 10 & 2  & 11 & 3  & 12 & 4  & 13 & 5  & 14 & 6  & 15 & 7  & 16 \\
10 & 2  & 11 & 3  & 12 & 4  & 13 & 5  & 14 & 6  & 15 & 7  & 16 & 8  & 0  & 9  & 1  \\
12 & 4  & 13 & 5  & 14 & 6  & 15 & 7  & 16 & 8  & 0  & 9  & 1  & 10 & 2  & 11 & 3  \\
14 & 6  & 15 & 7  & 16 & 8  & 0  & 9  & 1  & 10 & 2  & 11 & 3  & 12 & 4  & 13 & 5 \\
16 & 8  & 0  & 9  & 1  & 10 & 2  & 11 & 3  & 12 & 4  & 13 & 5  & 14 & 6  & 15 & 7 \\
1  & 10 & 2  & 11 & 3  & 12 & 4  & 13 & 5  & 14 & 6  & 15 & 7  & 16 & 8  & 0  & 9 \\
3  & 12 & 4  & 13 & 5  & 14 & 6  & 15 & 7  & 16 & 8  & 0  & 9  & 1  & 10 & 2  & 11 \\
5  & 14 & 6  & 15 & 7  & 16 & 8  & 0  & 9  & 1  & 10 & 2  & 11 & 3  & 12 & 4  & 13 \\
7  & 16 & 8  & 0  & 9  & 1  & 10 & 2  & 11 & 3  & 12 & 4  & 13 & 5  & 14 & 6  & 15 \\
9  & 1  & 10 & 2  & 11 & 3  & 12 & 4  & 13 & 5  & 14 & 6  & 15 & 7  & 16 & 8  & 0  \\
11 & 3  & 12 & 4  & 13 & 5  & 14 & 6  & 15 & 7  & 16 & 8  & 0  & 9  & 1  & 10 & 2  \\
13 & 5  & 14 & 6  & 15 & 7  & 16 & 8  & 0  & 9  & 1  & 10 & 2  & 11 & 3  & 12 & 4  \\
15 & 7  & 16 & 8  & 0  & 9  & 1  & 10 & 2  & 11 & 3  & 12 & 4  & 13 & 5  & 14 & 6 \\
\end{array} \right]\\
\end{array}
$$

\begin{center}
{\bf Figure IV}
\end{center}
\vspace{8pt}

\doublespace
Just as for the 2 and 3 dimensional cases, if $H =~ <\alpha_1,
\alpha_2, \alpha_3,\alpha_4>$ is a 4-D pandiagonal latin cube,
then so is $H^\prime$ with elements given by
$$  H^\prime_{ijk\ell} = H_{ijk\ell} + \alpha_1
    \equiv (\alpha_1 + 1) i + \alpha_2j + \alpha_3 k+ \alpha_4 \ell~~(mod~n)~.
$$
Again we can imagine ``tiling" all of 4-D space with $H$
and be assured that any $n\times n \times n \times n$ subcube in
this space will be a 4-D pandiagonal latin cube, with its elements
a permutation of the original elements of $H$.

\pagebreak

\noindent
{\small \bf III.  Magic, Pandiagonal Squares and Cubes in Two
and Three Dimensions}

With the results of the preceding
section available, it is now straightforward [{\bf 3}],[{\bf 6}]
to generate  {\it magic}
pandiagonal
cubes in three and four dimensions.
To set the stage we first illustrate the matter in
two dimensions.

\noindent
{\bf  A.  Magic Pandiagonal Squares}

Recall that a magic pandiagonal square of order $n$ is an
arrangement, without repeats, of the integers
$(0,1,2,...n^2-1)$ on the $n^2$ lattice points of an $n\times n$
array so that the sum of the elements in each of the
$n-$rows, $n-$columns, $n-$diagonals (including the $n-2$
broken diagonals) has the same value.
(As above let us use the generic {ROW} to represent any
one of these rows, columns, diagonals, etc. )
This common sum of the {ROWS} has been given in (1.1).

\noindent
{\bf Theorem 3.1}  [{\bf 9}, p. 178]  Let $L^{(1)} =~ <\alpha_1, \alpha_2>$,
$L^{(2)} = <\beta_1,\beta_2>$ be two orthogonal
pandiagonal latin squares that separately satisfy the conditions of
theorems (2.1) and (2.2) and define an $n\times n$ array
$M^{(2)}$ with elements $M_{ij}^{(2)}$ by
$$  M^{(2)} = nL^{(1)} + L^{(2)}
   \eqno{(3.1)}
$$
Then $M^{(2)}$ is a magic pandiagonal square.

\noindent
{\bf Proof:}  Since the elements of $L^{(1)}$ and $L^{(2)}$
range over $S = \{0,1,2,...n-1\}$ it follows from (3.1)
that each of the elements $M_{ij}^{(2)}$ must be
one of the integers $0,1,2,...n^2-1$.  Further, since the
sum of the elements in each {ROW} of the pandiagonal
$L^{(1)}$ and $L^{(2)}$, is $0+1+2+...+n-1$ $= n(n-1)/2$
it follows that the sum of each of the $4n$ ROWS of $M^{(2)}$ is
$$  n[n(n-1)/2] + n(n-1)/2 = n(n^2-1)/2 = \sigma_2
$$
with $\sigma_2$ defined in (1.1).
Finally since $L^{(1)}$ and
$L^{(2)}$ are orthogonal, no two elements of $M^{(2)}$
can be the same and the theorem is established{\bf .}

Since $L^{(1)}$ and $L^{(2)}$
are pandiagonal LS's for which we know $n \geq 5$, only
for $n \geq 5$ can $M^{(2)}$ in (3.1) be a magic pandiagonal square.
Further if we tile
the plane with $M^{(2)}$ any $n\times n$ subsquare will
also be a pandiagonal magic square since $L^{(1)}$ and $L^{(2)}$
have this same tiling property [{\bf 7}].  (See the discussion at the
end of IIA.)

Since the elements of the orthogonal pandiagonal latin squares
$L^{(1)}$ and $L^{(2)}$
may be permuted among themselves without changing their essential
properties, it follows that from a given pandiagonal LS we may generate $n$!
different versions.  This leads to the number $N_2$ of pandiagonal
magic squares obtainable by our method to be
$$  N_2 = \ell (n) (n!)^2
   \eqno{(3.2)}
$$
where $\ell(n)$ is a low order polynomial in $n$.  For
$n=5$, since the only independent LS's are
$<1,2>$ and $<1,3>$ (Theorem 2.3), and since we may interchange
their roles in (3.1)
$\ell(5) = 2$ and we obtain consistent with the results of
Rosser and Walker [11]
$$ N_2 = 2880
$$
Clearly because of the $(n!)^2$ factor
$N_2$ rises very rapidly
with $n$.  For $n=11$ for example, $N_2$ is $\sim$
$9.0 \times 10^{16}$.

\noindent
{\bf B.  Pandiagonal Magic Cubes}

A magic pandiagonal $n\times n\times n$ cube is an arrangement,
without repeats, of the integers $0,1,2,...n^3-1$ onto the
$n^3$ lattice sites of a cube so that the sums along each {ROW}
(i.e. along each of the $n^2$-rows, $n^2$-columns, $n^2$-files
and the $n(3n+6)$ diagonals (including the broken diagonals)
in each of its $(3n+6)$ squares
are the same.  The common sum along the {ROWS} is given by
$\sigma_3$ in (1.2).

\noindent
{\bf Theorem (3.2):}  Let $C^{(q)} =~ <\alpha_{1q},\alpha_{2q},
\alpha_{3q}>$, $q = 1,2,3$,
be three orthogonal pandiagonal $n\times n\times n$
latin cubes that satisfy the conditions of theorems (2.5),(2.6)
and define an $n\times n\times n$ array $M^{(3)}$ with elements
$M_{ijk}^{(3)}$ by
$$  M^{(3)} = n^2 C^{(1)} + nC^{(2)} + C^{(3)}.
   \eqno{(3.3)}
$$
Then $M^{(3)}$ is a magic pandiagonal cube.

\noindent
{\bf Proof:}  Since the elements of $C^{(1)}, C^{(2)}$ and $C^{(3)}$
range over
\linebreak
$S = \{0,1,2,...n-1\}$ it follows from (3.3) that each element
of $M^{(3)}$ must be one of the integers $\{0,1,2,...n^3-1\}$.
Consider now any ROW of
$M^{(3)}$.
The sum of the elements of this ROW is, according to (3.3), given by
$$  n^2[n(n-1)/2] + n[n(n-1)/2] + n(n-1)/2 = n(n^3-1)/2
$$
and yields the value $\sigma_3$ in (1.2).
Thus the sum of the elements
in each ROW of the cube is the same and since $C^{(1)}$,
$C^{(2)}$ and $C^{(3)}$
are orthogonal, it follows from (3.3) that no two
elements of $M^{(3)}$ can be the same.  The theorem is
thus established{\bf .}

According to the discussion in IIB, since the C's
are pandiagonal latin cubes, for which we know $n \geq 11$,
$M^{(3)}$ will exist only for $n \geq 11$.

As for the two dimensional case, if we allow for the interchange
in $M^{(3)}$,
%in $C^{(q)}$,
of $\alpha_{1q}, \alpha_{2q}, \alpha_{3q}$
$(q = 1,2,3)$ with each other, and of permuting the symbols
in each of $C^{(1)}$, $C^{(2)}$ and $C^{(3)}$ independently
we can conclude that $N_3$, the number of cubes $M^{(3)}$
in (3.3) is given by
$$  N_3 = \ell _3(n)(n!)^3
   \eqno{(3.4)}
$$
where $\ell_3$ is an appropriate polynomial in $n$.
For $n = 11,13$ and 17, $(n!)^3$ assumes the approximate values
$6.4 \times 10^{22}$, $2.4 \times 10^{29}$ and $4.5 \times 10^{43}$
respectively.  Because of this rapid rise of $N_3$ with $n$,
we can anticipate that the factor $\ell_3(n)$
will not affect this variation, qualitatively.

As for the analogous 2-D case, it follows that by virtue of
the ``tiling" properties of pandiagonal latin cubes (see IIB),
we can also tile 3-D space with any magic pandiagonal cube,
and be assured that any $n\times n\times n$ cube selected
out of this infinite array will also be a magic pandiagonal cube.
It will differ from the  original cube in that the elements
of its underlying latin cubes will have
been permuted.

\noindent
{\bf Four Dimensional Pandiagonal Magic Cubes}

By analogy to the above, we define a magic four dimensional
pandiagonal cube, as an arrangement, without repeats,
of the integers $0,1,2,...n^4-1$ among the $n^4$ lattice sites of
an $n\times n\times n \times n$ cube so that the sum of the elements in
each {ROW}
of the
4-D cube, has the same value $\sigma_4$ as given in (1.3).

\noindent
{\bf Theorem (3.3):}  Let $Q^{(q)} =~ <\alpha_{1q},\alpha_{2q},
\alpha_{3q}, \alpha_{4q}>$, $q = 1,2,3,4$ be four orthogonal
pandiagonal $n\times n\times n \times n$ latin cubes that satisfy
the conditions of theorems (2.6) and (2.7) and define an
$n\times n\times n\times n$ array $M^{(4)}$ with elements
$M_{ijk\ell}^{(4)}$ by
$$  M^{(4)} = n^3 Q^{(1)} + n^2Q^{(2)} + nQ^{(3)} + Q^{(4)}
  \eqno{(3.5)}
$$
Then $M^{(4)}$ is a magic, pandiagonal four dimensional cube.

\noindent
{\bf Proof:}  Since the elements of each of $Q^{(1)}, Q^{(2)}, Q^{(3)}$
and $Q^{(4)}$ range over
\linebreak
$S = \{0,1,2,...n-1\}$, it follows from (3.5)
that the elements of $M_{ijk\ell}^{(4)}$ assume values from the set
$\{0,1,2,...n^4-1\}$.
%Consider now any cube contained in $M^{(4)}$,
%and any square contained in this cube and any {ROW} contained in
%this square.
Further, since the $Q$'s are pandiagonal it follows that the
sum of the elements in any ROW of $M^{(4)}$ is given according to (3.5) by

$$  n^3[n(n-1)/2] + n^2[n(n-1)/2] + n[n(n-1)/2]
     + n(n-1)/2 = n(n^4-1)/2
$$
and this is $\sigma_4$ in (1.3).  Thus the sum of the elements
of any {ROW} of $M^{(4)}$ is given by (1.3).  Finally
because of the assumed orthogonality of the $Q$'s, it follows that
no two elements of
$M^{(4)}$
are the same and $M^{(4)}$
is a magic pandiagonal 4-D cube without repeats.
The theorem is established{\bf .}
%each of the $n^4$ numbers $0,1,2,....n-1$ appears at one
%of the $n^4$ lattice sites of $M^{(4)}$.

According to the discussion in IIC, since the $Q$'s are 4-D
pandiagonal latin cubes, for which we found
$n \geq 17$, it follows that $M^{(4)}$ will exist only for $n \geq 17$.

The number $N_4$ of different $M^{(4)}$'s in (3.5) can be
estimated as above to vary for large values of $n$ as
$$  N_4 \cong (n!)^4
$$
so that for $n=17, 19$ and 23, $N_4$ assumes the approximate
and rapidly growing values of $n$ of $1.6 \times 10^{58}$, $2.19 \times 10^{68}$ and
$4.5 \times 10^{89}$ respectively.

It is also possible, using the
tiling properties of the $Q's$
to
``tile" $M^{(4)}$ throughout four dimensional space and
obtain a pandiagonal magic 4-D cube by selecting any
$n\times n \times n \times n$ cube in this space.
Such a cube will be the same as the original cube but with
the elements of its underlying latin cubes
permuted.

Evidentally, these arguments are extendable to
dimensions higher than 4, but the resulting constraints on the $\alpha$
parameters,
analogous to those in (2.9) in three dimensions and those
in (2.14) in four, can be expected to become increasingly
complex.

\pagebreak

\begin{center}
{\bf REFERENCES}
\end{center}

\begin{enumerate}
%1
\item Brian Alspach and Katherine Heinrich, ``Perfect Magic Cubes
of Order 4m", The Fibonacci Quarterly, {\bf 19} (1981): 97-106.
%2
\item W.S. Andrews, ``Magic Squares and Cubes", New York;
Dover, 1960.
%3
\item W. Rouse Ball and H.M.S. Coxeter ``Mathematical Recreations Essays",
U. of Toronto Press, Cambridge, 1974; chapter VII.
%4
\item J. Denes and A.D. Keedwell, Latin Squares and their
Applications; New York Academic Press, 1974.
%5
\item Maurice Kraitchik, ``Mathematical Recreations", New York, Dover, 1942.
%6
\item Leonardi Euler, ``De quatratis magicis", Opera Omni,
S\'erie {\bf 1}, 7 (1923) 441-457; and
``Recherches sur un novelle esp\'ece de quarres magiques,"
Opera Omni, s\'erie {\bf 1}, 7, (1923), 291-392.
%7
\item Martin Gardner, ``Mathematical Games",
{\it Scientific American}, 234, 1 (January 1976), 118-123.
%8
\item Martin Gardner, ``Mathematical Games", {\it Scientific
American}, {\bf 234}, 2 (February 1976) 122-127.
%9
\item Charles F. Laywine and Gary L. Mullen,
``Discrete Mathematics Using Latin Squares",
New York, Wiley 1998.
%10
\item Paul. C. Pasles, ``The Lost Squares of Dr. Franklin":
Amer. Math. Monthly, {\bf 108}, (2001) 489-511.
%11
\item Barkeley Rosser and R.J. Walker, ``The Algebraic
Theory of Diabolic Magic Squares", Duke Math. J.,
5, (1939): 705-728.
%12
\item Bayard E. Wynne, ``Perfect Magic Cubes of Order Seven"
J. Recreational Mathematics, {\bf 8}:4, 1976, 285-293.
\end{enumerate}

\end{document}